\documentclass[11pt,a4paper]{article}
\usepackage[latin1]{inputenc}
\usepackage{indentfirst}
\usepackage{amsmath}
\usepackage{amsfonts}
\usepackage{amsthm}
\usepackage{amssymb}
\usepackage{graphics}
\usepackage{graphicx}

\def\H{\mathcal{H}}

\def\DD{\mathbb{D}}

\def\N{\mathbb{N}}

\def\P{\mathbb{P}}
\def\E{\mathbb{E}}

\def\D{\mathcal{D}}

\def\R{\mathcal{R}}
\def\Rp{\overset{\rightarrow }{\scriptstyle{ \mathcal{R}}}}
\def\Rpt{\overset{\longrightarrow }{\scriptstyle{ \mathcal{R}(t)}}}
\def\Rpl{\overset{\longrightarrow }{\scriptstyle{ \mathcal{R}(\alpha,1)}}}
\def\Rpa{\overset{\longrightarrow }{\scriptstyle{ \mathcal{R}(\alpha,\lambda)}}}

\def\Rm{\overset{\leftarrow }{\scriptstyle{\mathcal{R}}}}
\def\Rmt{\overset{\longleftarrow }{\scriptstyle{ \mathcal{R}(t)}}}
\def\Rml{\overset{\longleftarrow }{\scriptstyle{ \mathcal{R}(\alpha,1)}}}
\def\Rma{\overset{\longleftarrow }{\scriptstyle{ \mathcal{R}(\alpha,\lambda)}}}

\def\C{\mathcal{C}}
\def\RRR{\mathbb{R}}

\def\t{\textrm}
\def\w{\widetilde}

\def\B{\B }
\def\taup{\overset{\rightarrow }{\scriptstyle{\tau}}}
\def\taum{\overset{\leftarrow}{\scriptstyle{\tau}}}
\def\B{\mathbf{B} _0}
\def\ind{{\mathchoice {\rm 1\mskip-4mu l} {\rm 1\mskip-4mu l}
{\rm 1\mskip-4.5mu l} {\rm 1\mskip-5mu l}}}

\newcommand{\be} {\begin{equation}}
\newcommand{\ee} {\end{equation}}
\newcommand{\bea} {\begin{eqnarray}}
\newcommand{\eea} {\end{eqnarray}}
\newcommand{\Bea} {\begin{eqnarray*}}
\newcommand{\Eea} {\end{eqnarray*}}

\newtheorem{Thm}{Theorem}
\newtheorem{Lem}{Lemma}

\newtheorem{Pte}{Proposition}
\newtheorem{Cor}{Corollary}
\theoremstyle{definition} \newtheorem*{key}{Key words}
\theoremstyle{definition} \newtheorem*{ams}{A.M.S. Classification}
\theoremstyle{definition} \newtheorem{ex}{Example}
\theoremstyle{definition} 

\theoremstyle{remark}\newtheorem{Rque}{Remark}
\theoremstyle{remark}\newtheorem{fig}{Figure}

\begin{document}
\title{On a model for the storage of files on a hardware I :\\
Statistics at a fixed time and asymptotics.}
\author{Vincent Bansaye \footnote{Laboratoire de Probabilités et Modèles Aléatoires. Université Pierre et Marie Curie et C.N.R.S. UMR 7599. 175, rue du
Chevaleret, 75 013 Paris, France.
$\newline$
\emph{e-mail} : bansaye@ccr.jussieu.fr }}
\maketitle
\vspace{3cm}
\begin{abstract}We consider  a generalized version in continuous time of the parking problem
of Knuth. Files arrive following a Poisson point process and are stored on a hardware identified with the real line. We
specify the
distribution of the space of unoccupied locations  at a fixed time and give its asymptotics
when the hardware is becoming full.
\end{abstract}
\begin{key}Parking problem. Data storage. Random covering. Lévy processes. Regenerative sets. Queuing theory.
\end{key}
\begin{ams} 60D05, 60G51, 68B15.
\end{ams}

\section{Introduction}

We consider  a generalized version in continuous time of the
original parking problem of Knuth. Knuth was interested by the
storage of data on a  hardware represented by a circle with  $n$
spots. Files arrive successively  at locations chosen uniformly
among these $n$ spots. They are stored in the first free spot at
the right of their arrival point (at their arrival point if it is
free). Initially Knuth worked on  the hashing of data
(see e.g. \cite{chass, Flaj, Foa})  : he
 studied the distance  between the spots where the files arrive and the spots where they are stored. Later
Chassaing and Louchard
\cite{chl} have described the evolution of the
largest block of data in such coverings when $n$ tends to infinity. They observed a phase transition at the stage
where the hardware is almost full, which is
 related to the additive coalescent.
Bertoin and Miermont \cite{bem} have extended these results to files of random sizes which arrive uniformly on the circle. \\

We consider here a continuous time version of this  model where  the
hardware is large  and  now identified  with the real line. A file labelled $i$
of length (or size) $l_i$  arrives at time $t_i\geq 0$  at location
$x_i\in \RRR$.
 The storage of this file uses the free portion of size $l_i$ of the real line at the right of
$x_i$ as close to $x_i$ as possible (see Figure \ref{fig1}). That
is, it covers $[x_i,x_i+l_i[$ if this interval  is free at time
$t_i$. Otherwise it is shifted to the right until a free space is
found and it may be split into several parts which are stored in the
closest free portions. We require absence of memory for the storage of  files, uniformity of the location
where they arrive and
identical distribution of the sizes. Thus, we model the arrival of
files by a Poisson point process (PPP) :
 $\{(t_i,x_i,l_i)\  : \ i\in \N\}$ is a  PPP  with intensity $\t{d}t\otimes \t{d}x\otimes \nu(\t{d}l)$ on
$\RRR^+\times \RRR\times \RRR^+$. We  denote $\t{m}:=\int_0^\infty l \nu(\t{d}l)$ and assume $\t{m} < \infty$. So $\t{m}$
is the mean of the total sizes of files which arrive during a unit interval time on some interval with unit length.  \\

We begin by constructing this random covering (Section 2.1). The first questions which
arise and are treated here concern  statistics at a fixed time.
What is the distribution of the covering at a fixed time ? At what time  the hardware becomes full ?  What are
the asymptotics of the covering at this saturation time ? What is the length of the largest block on a part of the hardware ? \\

It is quite easy to see that the hardware becomes full at a
deterministic time equal to $1/\t{m}$. In Section 3, we give
some geometric properties of the covering and characterize the
distribution of the covering at a fixed time  by
giving the joint distribution of the block of data strad\t{d}ling $0$
and the free spaces on both sides of this block. The results given in this section
will be  useful for
the problem of the dynamic of the covering considered in
$\cite{vbb}$, where we characterize the evolution in time of a typical data block. Moreover, using  this characterization,  we
determine the asymptotics  of the covering $\C(t)$ at the saturation time
$1/\t{m}$ (Theorem \ref{asympt}). \\
 By the same method,  we determine the 
asymptotic regime of the hardware restricted to $[0,x]$ rescaled to $[0,1]$ at saturation time (Theorem \ref{thtrst}). We derive then the  asymptotic of the 
largest block
of the hardware restricted to $[0,x]$ when $x$ tends to infinity. As expected, we recover the 
phase transition observed by Chassaing and Louchard in \cite{chl}.\\

As we look at $\C(t)$ at a fixed $t$, it does not depend on the order of of files before 
time $t$. Thus if  $\nu$ is finite, we can view the files which arrive before time $t$ 
as customers :  the size of the file $l$ is the service time and the location $x$ where 
the file arrives is the arrival time of the customer. We are then in the framework of $M/G/1$ 
queuing model in the stationary regime and the covering $\C(t)$ is the union of busy periods 
(see e.g. Chap 3 in \cite{Prabhu}). Thus, results of Section 3 for finite $\nu$  follow easily from known results on $M/G/1$.  For infinite $\nu$, results are the same but busy cycle is not defined and proofs are different and proving asymptotics on random sets  need results  about Lévy processes.  Moreover, as far as we know, the longest
busy period and more generally asymptotic regimes on 
 $[0,x]$ when $x$ tends to infinity and $t$ tends to the saturation time has not been considered in queuing model. 
$\newline$ 

\section{Preliminaries}

Throughout this paper, we use the classical notation  $\delta_x$ for the Dirac mass at $x$  and
$\N=\{1,2,..\}$.  If
$\R$ is a measurable subset of $\RRR$,  we denote by $\mid \R\mid $ its Lebesgue measure and
 by $\R^{cl}$ its closure. For every $x\in\RRR$, we denote by $\R-x$ the set $\{y-x  : x \in \R\}$ and  \\
\be \label{defgd} g_x(\R)=\t{sup}\{y\leq x : y \in \R\}, \ \ \ \ \
d_x(\R)=\t{inf}\{y>x : y \in \R\}. \ee
If $I$ is a closed interval of $\RRR$, we denote by $\H(I)$ the
space of closed subset of $I$. We endow $\H(I)$ with the Hausdorff
distance $d_H$  defined  for all $A,B \subset \RRR$ by :
$$d_H(A,B)=\t{max}\big(\sup_{x\in A} d(x,B),     \sup_{x\in B} d(x,A)\big), \quad \t{where} \
d(x,A)=\t{inf}\{1-e^{-\vert x-y\vert} :  y \in A \}.$$
 The topology induced by this distance is the topology of Matheron \cite{mat} : a sequence $\R_n$ in $\H(I)$
converges to $\R$ iff for each open set $G$ and each compact $K$,
\bea
&&\R\cap G\ne \varnothing \ \ \  \t{implies} \ \ \ \R_n\cap G\ne\varnothing \ \ \ \t{for} \  n \ \t{large enough} \nonumber \\
&&\R\cap K= \varnothing \ \ \ \t{implies} \ \ \   \R_n\cap K=\varnothing \ \ \ \t{for} \  n  \ \t{large enough} \nonumber
\eea
It is also the topology induced by the Hausdorff metric on  a compact using arctan($\R\cup\{-\infty,\infty\}$)
 or the Skorokhod metric using the class of 'descending saw-tooth functions' (see \cite{mat} and \cite{frm} for details). \\

\subsection{Construction of the covering $\C(t)$}

First, we present a deterministic construction of the covering $\C$ associated with a given sequence of files labeled by $i \in \N$. The file
labeled by  $i\in\N$ has size $l_i$ and arrives after the files labeled by $j\leq i-1$, at location $x_i$ on the real line. Files
are stored following the process described in the Introduction and $\C$ is the portion of line which is used for the storage.  We
begin by constructing
the covering $\C^{(n)}$  obtained by considering only the first $n$ files, so that $\C$  is obtained as  the increasing union of these
coverings.
A short though (see Remark \ref{commut}) enables us to see that the covering $\C$  does not depend on the order of arrival of the
files. This  construction of $\C$ will then be applied to  the construction of our random covering at a fixed time $\C(t)$ by
considering files arrived before time $t$.
$\newline$
\begin{fig}
\label{fig1}
Arrival and storage of the $5$-file  and representation of $Y^{(5)}$. The first four files have been stored without
splitting and are represented by the black rectangles.
$\includegraphics[scale=0.4]{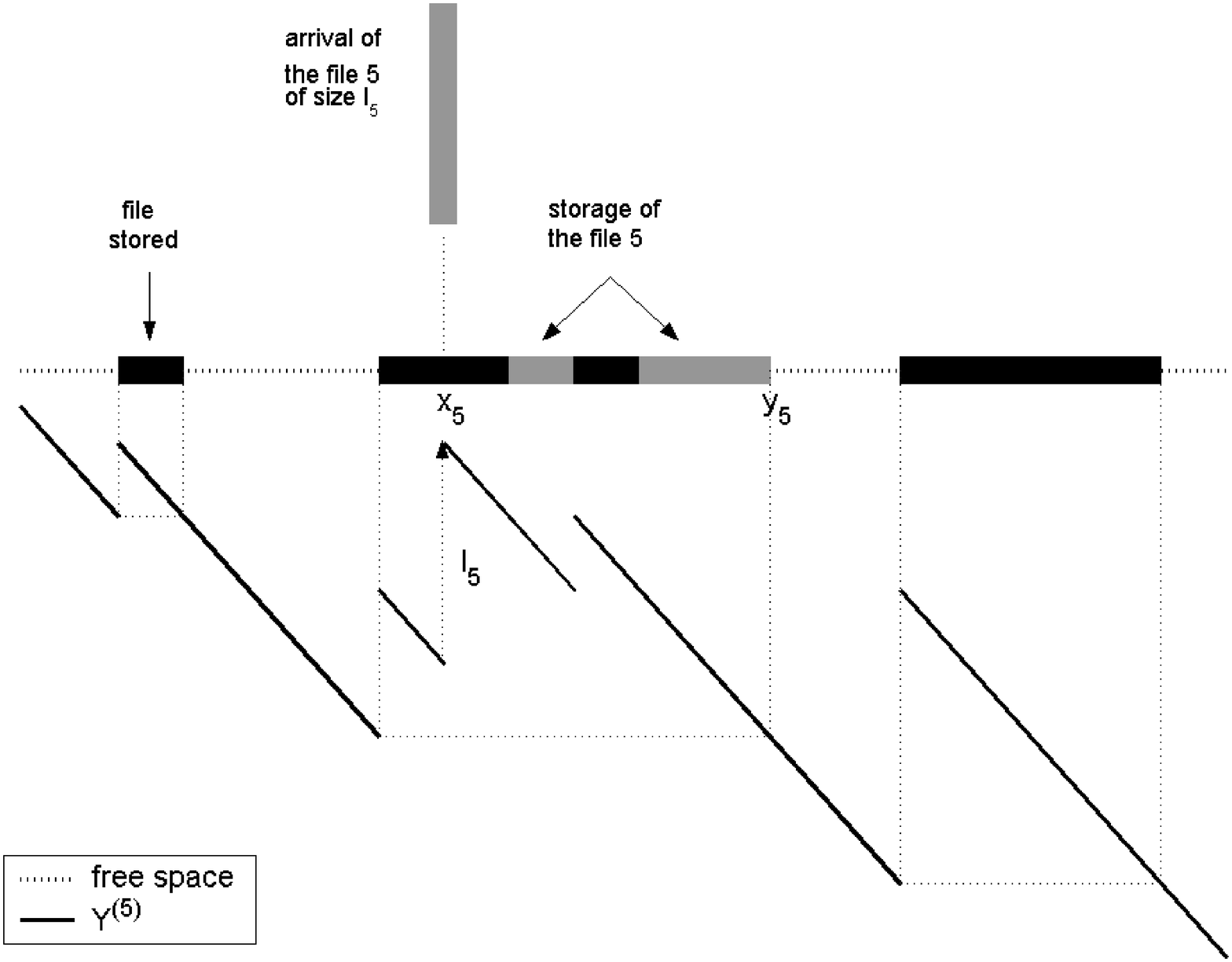}$
\end{fig}

We define $\C^{(n)}$ by induction  :
$$\C^{(0)} :=\varnothing \ \ \ \ \ \  \C^{(n+1)}:=\C^{(n)}\cup [x_{n+1},y_{n+1}[ $$
where $y_{n+1}=\t{inf}\{y\geq 0, \ \mid \R^{(n)}\cap [x_{n+1},y[\mid =l_{n+1}\}$ and $\R^{(n)}$ is the
 complementary set of $\C^{(n)}$ (i.e. the free space of the real line). So $y_{n+1}$ is the right-most
point which is used for storing the $(n+1)$-th file.\\ \\
Now we consider the quantity of data over $x$,  $R^{(n)}_x$, as the quantity of
data which we have tried to store at the location $x$
(successfully or not) when $n$ files are stored.  These data are
the data fallen in $[g_x(\R^{(n)}),x]$ which could not be stored
in $[g_x(\R^{(n)}),x]$, so  $R^{(n)}_x$ is defined by
$$R^{(n)}_x:= \sum_{\substack{i\leq n \\ x_i \in [g_x(\R^{(n)}),x]}} l_i \ \ - \ \ (x-g_x(\R^{(n)})).  $$
This  quantity can be expressed using  the function $Y^{(n)}$, which  sums the sizes of the files arrived
at the left of a point $x$ minus   the drift term $x$.  It is thus defined by $Y^{(n)}_0=0$ and
$$Y^{(n)}_b-Y^{(n)}_a= \sum_{\substack{i\leq n \\  x_i \in ]a,b]}} l_{i} \ \ -  \ \ (b-a) \qquad \t{for} \ a<b.$$
Introducing also its infimum function  defined for $x\in\RRR$ by  \ \ $I^{(n)}_x:=\t{inf}\{Y^{(n)}_y  : y\leq x\}$, we get
the following expression.
$\newline$
\begin{Lem}
\label{Re}
For every $n\geq 1$, we have  $R^{(n)}=Y^{(n)}-I^{(n)}$.
\end{Lem}
\begin{proof}
Let $x\in\RRR$. For every $y\leq  x$, the quantity of data over $x$ is at least the quantity of data fallen  in $[y,x]$ minus $y-x$, i.e.
$$R^{(n)}_x\geq \sum_{\substack{i\leq n \\ \ x_i \in [y,x]}} l_i \  -(x-y)$$
and by definition of $R_x^{(n)}$, we get :
$$
 R^{(n)}_x= \t{sup}\{ \sum_{\substack{i\leq n \\ x_i \in [y,x]}} l_i \ - \ (x-y) : \ y\leq x\}=  \t{sup}\{Y^{(n)}_x-Y^{(n)}_y : \ y\leq x\}.$$
Then $R^{(n)}_x=Y^{(n)}_x-I^{(n)}_x $.
\end{proof}
$\newline$
As a consequence, the covered set when the first $n$ files are stored is given by
\be
\label{cover}
\C^{(n)}=\{ x\in \RRR : \ Y^{(n)}-I^{(n)}>0\}.
\ee

We are now able to investigate the situation when $n$ tends to infinity under the following   mild condition
\be
\label{condH}
\qquad \qquad \forall \ L\geq 0, \ \ \sum_{x_i\in [-L,L]} l_i <\infty,
\ee
which means that  the quantity of data arriving on a compact set is finite. We introduce  the   function $Y$ defined on $\RRR$ by $Y_0=0$ and
$$Y_b-Y_a= \sum_{  x_i \in ]a,b]} l_{i} \ \ -  \ \ (b-a) \qquad \t{for} \ a<b$$
and its infimum $I$ defined for $x\in\RRR$ by \ \ $I_x:=\t{inf}\{Y_y : y\leq x\}$. \\ \\
As expected, the covering $\C:=\cup _{n\in \N} \C^{(n)}=\cup _{n\in \N} \{ x\in \RRR : \ Y^{(n)}-I^{(n)}>0\}$ is given by
\begin{Pte}
\label{disc}
- If $\lim_{x\rightarrow -\infty} Y_x=+\infty$, then  $\C=\{ x\in \RRR : \ Y_x-I_x>0\}\ne \RRR$.\\ \\
- If $\liminf_{x\rightarrow -\infty} Y_x=-\infty$, then  $\C=\{ x\in \RRR : \ Y_x-I_x>0\}= \RRR$.\\
\end{Pte}
\begin{Rque}
\label{commut}
This result ensures that the covering does not depend on the order of arrival of files.
\end{Rque}
\begin{proof}
Condition (\ref{condH}) ensures that $Y^{(n)}$ converges  to $Y$ uniformly on every compact set of $\RRR$. \\ \\
$\bullet$ \ If $\lim_{x\rightarrow -\infty} Y_x=+\infty$, then for every  $L \geq 0$, there exists $L'\geq L$ such
that $I_{-L'}=Y_{-L'}$. Moreover $Y_x\leq Y^{(n)}_x$ if $x\leq 0$. So :
$$Y^{(n)}_{-L'}\stackrel{n \rightarrow \infty}{\longrightarrow }Y_{-L'}=I_{-L'} \ \ \ \t{and} \ \ \ \ I_{-L'}\leq I^{(n)}_{-L'}\leq Y^{(n)}_{-L'} $$
Then $I^{(n)}_{-L'}\stackrel{n \rightarrow \infty}{\longrightarrow }I_{-L'}$. As $Y^{(n)}$ converges  to $Y$ uniformly on
$[-L',L']$, this  entails that for every $x$ in $[-L,L]$,  $\t{inf}\{Y^{(n)}_y, -L'\leq y\leq x\}
\stackrel{n \rightarrow \infty}{\longrightarrow }\t{inf}\{Y_y, -L'\leq y\leq x\}$. Then,
$$ \ I^{(n)}_{x}=I^{(n)}_{-L'}\wedge \t{inf}\{Y^{(n)}_y, -L'\leq y\leq x\} \stackrel{n \rightarrow \infty}{\longrightarrow }
I_{-L'}\wedge \t{inf}\{Y_y, -L'\leq y\leq x\} =I_{x}.$$
So $Y^{(n)}_x-I^{(n)}_x \stackrel{n \rightarrow \infty}{\longrightarrow }Y_x-I_x$ and $Y^{(n)}_x-I^{(n)}_x$ increases when $n$
increases since it is equal to $R_x^{(n)}$,
the quantity of data over $x$  (see Lemma \ref{Re}). We conclude that there is the identity
$$\{ x\in \RRR, \ Y_x-I_x>0\}=\cup _{n\in \N}\{ x\in \RRR, \ Y^{(n)}_x-I^{(n)}_x>0\}=\C.$$ \\

Moreover $-L' \notin \{ x\in \RRR : \ Y_x-I_x>0\}$, so $\C=\{ x\in \RRR, \ Y_x-I_x>0\} \ne \RRR$. \\ \\
$\bullet$ If $\liminf_{x\rightarrow -\infty} Y_x=-\infty$, then for every $x \in \RRR$,
$$I_x=-\infty  \ \ \ \ \ \ \t{and} \ \ \ \ \ \  I^{(n)}_{x}\stackrel{n \rightarrow \infty}{\longrightarrow }-\infty$$
The first identity entails that  $\{ x\in \RRR, \ Y_x-I_x>0\}= \RRR$. As  $(Y^{(n)}_x)_{n \in \N}$ is bounded, the second
one implies that
there exists $n$ in $\N$ such that $Y^{(n)}_x-I^{(n)}_{x}>0$. Then we have also
$\cup _{n\in \N} \{ x\in \RRR : \ Y^{(n)}-I^{(n)}>0\}=\RRR$, which gives the result.
\end{proof}
$\newline$

Finally, we can  construct  the random covering associated with  a PPP. As the order of arrival of files has no importance, the
random covering $\C(t)$ at time $t$ described in Introduction is
obtained by the deterministic construction above by taking the subfamily of files $i$ which verifies $t_i\leq t$.
$\newline$

When files arrive according to a PPP,   $(Y_x)_{x\geq 0}$ is a Lévy process,  and
we recall now some results about Lévy processes and their fluctuations which will be useful in the rest of this work.
$\newline$

\subsection{Background on Lévy processes}

The results given in this section can be found  in  the Chapters VI and VII in \cite{lev} (there, statements are
made in terms
of the dual process $-Y$). We recall that a Lévy process is cà\t{d}làg process starting from $0$ which has
iid increments. A subordinator is an increasing Lévy process.

We consider in this section  a Lévy process $(X_x)_{x \geq 0}$ which has no negative jumps (spectrally positive Lévy process).  We
 denote by $\Psi$ its Laplace exponent which verifies for every $\rho\geq 0$ :
\be
\label{lapdef}
\E(\t{exp}(-\rho X_x))=\t{exp}(-x \Psi (\rho )).
\ee
We stress that this is not  the classical choice for the sign of the Laplace exponent of  Lévy processes with no negative
jumps and a negative drift such  as the process ($Y_x)_{x\geq 0}$ introduced in the previous section. However it is the classical
choice for   subordinators, which we will need. It is then convenient  to use this same definition for all
 Lévy processes which appear in this text. \\ \\

First, we consider the  case when  $(X_x)_{x \geq 0}$ has bounded variations. That is,
$$ X_x:=dx + \sum_{x_i \leq x} l_i, $$
where $\{ (x_i,l_i) : i \in \N \} $ is a PPP on $[0,\infty[\times [0,\infty]$
with intensity measure $ \t{d}x \otimes  \nu $ such that $ \int_0^{\infty}x \nu (\t{d}x) < \infty $. We call $\nu$  the Lévy
measure and   $d \in \RRR$  the drift.  Note that  $(Y_x)_{x\geq 0}$ is a subordinator
iff $d \geq 0$.\\ \\
Writing $\bar \nu$ for the tail of the measure $\nu$, the Lévy-Khintchine formula gives
\bea
\label{Lap1}
\Psi (\rho)&=& d \rho  +\int_0^{\infty} (1-e^{-\rho x}) \nu (\t{d}x), \\
\label{Lap2}
\frac{\Psi (\rho)}{\rho}&=&d+\int_0^{\infty} e^{-\rho x} \bar \nu  (x)\t{d}x, \\
\label{lim2}
\Psi'(0)&=&d+\int_0^{\infty} x \nu  (\t{d}x), \ \ \ \ \ \ \ \ \\
\label{lim1}
\lim_{\rho \rightarrow \infty} \frac{\Psi (\rho)}{\rho} =d   &\t{and}&   \lim_{\rho \rightarrow \infty} 
(\Psi (\rho)-d \rho) = \bar \nu (0).
\eea

$\newline$

Second, we consider the case when $\Psi$ has a right  derivative at $0$ with
\be
\label{condHH}
\Psi'(0)<0
\ee
meaning that $\E(X_1)<0$. And we consider  the infimum process  which has  continuous path  and the first
 passage time  defined for $x\geq 0$ by
 $$I_x=\t{inf}\{ X_y : 0\leq y\leq x \} \ \ \ ; \ \ \ \ \tau_x=\t{inf}\{z\geq 0 : X_z<-x \}.$$
 As $-\Psi$ is strictly convex and $-\Psi'(0)>0$, $-\Psi$ is strictly increasing from $[0,\infty[$  to $[0,\infty[$ and so
 is strictly positive on $]0,\infty]$. We write $\kappa : [0,\infty[ \rightarrow \mathbb{R} $ for
the inverse function of $-\Psi$ and we have (see \cite{lev}  Theorem 1 on page 189 and Corollary 3 on page 190) :
$\newline$
\begin{Thm}
\label{tau}
$(\tau_x)_{x\geq 0}$ is a subordinator with Laplace exponent $\kappa$.\\

Moreover the following identity  holds between measures  on $[0,\infty[\times [0,\infty[$  :
\be
\label{egalintro}
x\P(\tau_l \in \t{d}x)\t{d}l=l\P(-X_x \in \t{d}l)\t{d}x.
\ee
\end{Thm}
$\newline$
Note that if  $(X_x)_{x \geq 0}$ has  bounded variations, using (\ref{lim1}), we can write
\be
\label{kappa}
\forall \ \rho\geq 0, \quad \kappa(\rho)=-\frac{\rho}{d} + \int_0^{\infty} (1-e^{-\rho z}) \Pi (\t{d}z),
\ee
where $\Pi$ is a measure on $\RRR^+$ verifying (use ($\ref{lim1}$)  and  Wald's identity or  ($\ref{lim2}$)) :
\be
\label{intPi}
\bar \Pi(0)=-\frac{\bar \nu(0)}{d},  \qquad \int_0^{\infty} x\Pi(\t{d}x) =\frac{1}{d}-\frac{1}{d+\int_0^{\infty} x\nu(\t{d}x)}.
\ee
$\newline$
Now we introduce the supremum process defined for $x\geq 0$ by
$$S_x:=\t{sup}\{ X_y : 0\leq y\leq x \},$$
 and the a.s unique instant   at which
$X$ reaches this supremum on $[0,x]$ :
$$\gamma_x:=\t{inf}\{y \in [0,x] : X_y=S_x\}.$$
By duality, we have $(S_x,\gamma_x) \stackrel{d}{=} (X_x-I_x,x-g_x)$ where $g_x$ denotes the a.s unique instant at which $(X_{x^-})_{x\geq 0}$ reaches its overall infimum on $[0,x]$
(see Proposition 3 in \cite{lev} or \cite{mono} on page 25).
If  $T$ is an  exponentially distributed random time with parameter $q>0$ which is independent of $X$ and
$\lambda,\mu > 0$, then we have   (use \cite{lev} Theorem 5 on page 160 and Theorem 4 on page 191) :
\bea
 \E\big( \t{exp} (-\mu S_{T}- \lambda \gamma_{T})\big)&=&\frac{q(\kappa(\lambda+q)-\mu)}{\kappa(q)(q+\lambda+\Psi(\mu))} \nonumber \\
 &=&\t{exp}\big(\int_0^{\infty} \t{d}x \int_0^{\infty}\P(Y_x\in \t{d}y)(e^{-\lambda x -\mu y}-1)x^{-1}e^{-qx}\big) \nonumber
\eea
which gives
\bea
\label{sup}
\E\big( \t{exp} (-\mu S_{\infty}-\lambda \gamma_{\infty})\big)&=&\frac{1}{\kappa'(0)} \frac{\kappa(\lambda)-\mu}{\lambda+\Psi(\mu)} =-\Psi'(0) \frac{\kappa(\lambda)-\mu}{\lambda+\Psi(\mu)} \\
\label{sup2}
\E\big( \t{exp} (-\mu S_{\infty})\big)&=& \mu\frac{\Psi'(0)}{\Psi(\mu)}  \\
\label{fluc}
\E\big(\t{exp} (-\lambda \gamma_{\infty})\big)&=&\t{exp}\big(\int_0^{\infty} (e^{-\lambda x}-1)x^{-1}\P(X_x> 0)\t{d}x\big)
\eea

$\newline$

\section{Properties of the covering at a fixed time}

\subsection{Statement of the results}
Our purpose in this section is to specify the distribution of  the
covering $\C(t)$ and we will use the characterization of Section 2.1 and
results of Section 2.2.  In that view, following the previous section,
we consider  the process  $(Y^{(t)}_x)_{x\in \RRR}$   associated to the PPP $\{(t_i,l_i,x_i), i\in
\N\}$ and defined by
$$Y^{(t)}_0:=0 \quad ; \quad  Y^{(t)}_b-Y^{(t)}_a= \sum_{\substack{t_i\leq t \\ 
x_i \in ]a,b]}} l_{i} \ \ -  \ \ (b-a) \quad \t{for} \ a<b,$$
which has independent and stationary increments, no negative jumps and
bounded variation. Introducing also its infimum process defined for $x\in\RRR$ by
$$I^{(t)}_x:=\t{inf}\{Y^{(t)}_y : y\leq x\},$$
 we can give now
 a handy expression for the covering at a fixed time and obtain that the hardware becomes full at a deterministic
time equal to $1/\t{m}$ (see below for the proofs).
$\newline$
\begin{Pte}
\label{dico}
 For every $t<\emph{1/m}$, we have $\C (t)= \{ x\in \RRR : \ Y^{(t)}_x > I^{(t)}_x \} \ne \RRR$ a.s. \\
 For every $t\geq \emph{1/m}$, we have $\C (t)=\RRR$ a.s.
\end{Pte}

To specify the distribution of $C(t)$, it is equivalent and  more convenient to describe its complementary
set, denoted by $\R(t)$, which corresponds to the free
space of the hardware. By the previous proposition, there is the
identity  :
\be \label{comp} \R(t)=\{x \in \RRR :
Y^{(t)}_x=I^{(t)}_x\}. \ee
We begin by giving  some classical geometric properties  which will be  useful.
$\newline$
\begin{Pte}
\label{prob} For every  $t \geq 0$,  $\R(t)$  is  stationary, its  closure is symmetric in distribution  and it enjoys
the regeneration property :
For every $x \in \RRR$,  $(\R(t)-d_x(\R(t)))\cap [0,\infty[$ is independent of $\R(t)\cap ]-\infty,x]$ and
is distributed as  $(\R-d_0(\R(t)))\cap [0,\infty[$. \\
Moreover for every  $x \in \RRR$, $\P(x\in \C (t))=\emph{min}(1,\emph{m} t)$.
\end{Pte}
\begin{Rque}
Even though  the distribution of $\R(t)^{cl}$ is symmetric, the processes  $(\R(t)^{cl} : t \in [0,1/\t{m}])$ and
$(-\R(t)^{cl} : t \in [0,1/\t{m}])$  are quite different. For example, we shall  observe in \cite{vbb}
that the left extremity of the data block strad\t{d}ling $0$  is a Markov process but the right extremity  is not.
\end{Rque}
$\newline$

We want now to characterize the distribution of the free space $\R(t)$. For this purpose, we  need some notation.
 The drift of the  Lévy process  $(Y^{(t)}_x)_{x\geq 0}$
is equal to  $-1$, its  Lévy measure is equal to $t\nu$ and its
Laplace exponent $\Psi^{(t)}$ is then given by (see
($\ref{Lap1}$)) \be \label{Psi}
\Psi^{(t)}(\rho):=-\rho+\int_0^{\infty} \big( 1-e^{-\rho x} \big) t
\nu(\t{d}x). \ee
For sake of simplicity, we write,
recalling ($\ref{defgd}$),
$$ g(t):=g_0(\R(t)),    \qquad      d(t)=d_0(\R(t)),  \qquad l(t)=d(t)-g(t),$$
 which are respectively the left extremity, the right extremity and
the length  of the data block
strad\t{d}ling $0$, $\B (t)$.  Note that $g(t)=d(t)=0$ if $\B
(t)=\varnothing$. \\
We  work with $\R$
subsets of $\RRR$  of the form  $\sqcup_{n\in\N}[a_n,b_n[$ and we denote
by  $\w{\R}:=\sqcup_{n\in\N}[-b_n,-a_n[$ the symmetrical of
$\R$ with respect to $0$  closed at the left, open at the right. We consider the positive part (resp. negative part)
of $\R$ defined  by \bea
\label{defrp}
&&\Rp :=(\R-d_0(\R))\cap [0,\infty]=\bigsqcup_{n \in \N  : \ a_n\geq d_0(\R)}[a_n-d_0(\R),b_n-d_0(\R)[, \nonumber \\
\label{defrm}
&&\Rm := \overset{\rightarrow }{\w{\scriptstyle{\mathcal{R}}}}=\bigsqcup_{n \in\N :\ b_n\leq g_0(\R)} [g_0(\R)-b_n,g_0(\R)-a_n[. \nonumber
\eea
\begin{ex}
For a given $\R$ represented by the dotted lines, we give below $\Rp$ and $\Rm$, which are also represented by dotted lines. Moreover the endpoints of the data blocks containing $0$ are denoted by $g_0$ and $d_0$.

$\includegraphics[scale=0.4]{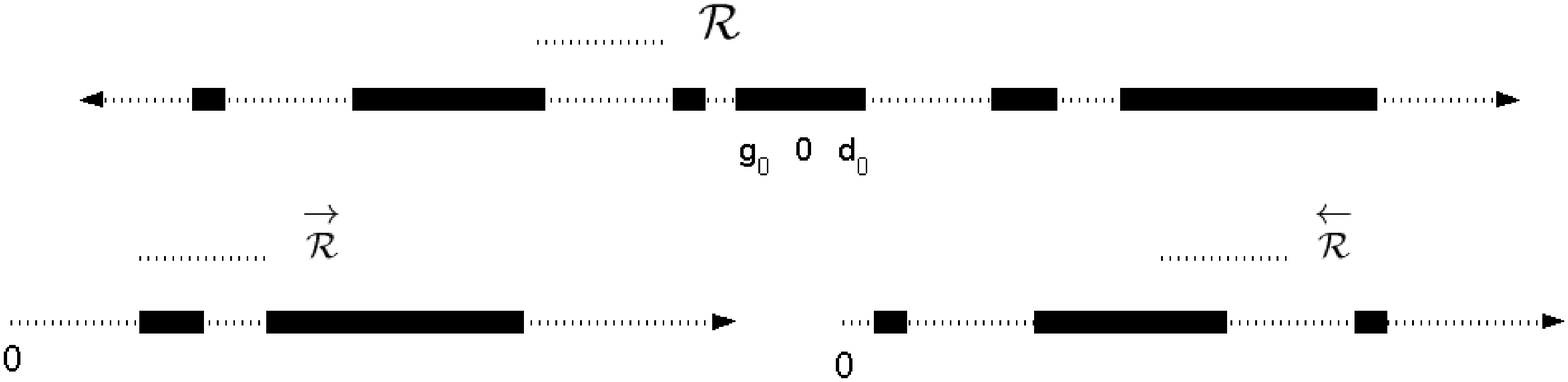}$
\end{ex}

Thus $\Rpt$ (resp. $\Rmt$) is the free space at the
right of $\B(t)$ (resp. at the left of  $\B(t)$, turned over,
closed at the left and open at the right). We have then the
identity \be \label{reprR}
\R(t)=(d(t)+\Rpt)\sqcup(\w{-g(t)+\Rmt}). \ee
Introducing also the processes  $(\taup^{(t)}_x)_{x\geq 0}$ and $(\taum^{(t)}_x)_{x\geq 0}$ defined by
$$\taup^{(t)}_x:=\t{inf}\{y\geq 0 : \vert\Rpt\cap [0,y]\vert>x\}, \qquad \taum^{(t)}_x:=\t{inf}\{y\geq 0 : \vert\Rmt\cap [0,y]\vert>x\},$$
enables us to describe  $\R(t)$ in the following way :
$\newline$
\begin{Thm}
\label{loi}
(i) The random  sets  $\Rpt$ and  $\Rmt$ are independent,  identically distributed and  independent of  $(g(t),d(t))$.\\
(ii) $\Rpt$ and  $\Rmt$ are   the range  of the subordinators $\taup^{(t)}$ and $\taum^{(t)}$  respectively
whose Laplace exponent $\kappa^{(t)}$ is the inverse function of $-\Psi^{(t)}$.\\
(iii) The  distribution of $(g(t),d(t))$ is specified by :
$$(g(t),d(t))=(-Ul(t),(1-U)l(t)),$$
$$\P(l(t) \in \emph{d}x)=(1-\emph{m}t)\big(\delta _0 (\emph{d}x)+\ind _{\{x>0\}} x\Pi^{(t)}(\emph{d}x) \big)$$
where $U$ uniform random variable on $[0,1]$ independent of $l(t)$ and
 $\Pi^{(t)}$ is the Lévy measure of $\kappa^{(t)}$.
\end{Thm}
$\newline$

We can then estimate the number of data blocks on the hardware. If $\nu$ has a finite mass, we write $N^{(t)}_x$ the  number of data blocks  of the hardware  restricted to $[-x,x]$ at time $t$. This quantity
 has a  deterministic asymptotic as $x$ tends to infinity  which is maximum at time $1/(2\t{m})$. In this sense, the number of blocks of the hardware
reaches a.s. its maximal at time $1/(2\t{m})$. More precisely,
\begin{Cor}\label{rq}If $\bar{\nu}(0)<\infty$, then for every $t\in[0,1/\emph{m}[$,
$$\lim_{x\rightarrow \infty}\frac{N^{(t)}_x}{2x}=\bar{\nu}(0)t(1-\emph{m}t) \quad \t{a.s.}$$
\end{Cor}
$\newline$

Moreover, we can describe here the hashing of data. We recall that
a file labelled by $i$ is stored at location $x_i$. In the hashing
problem,  one is interested by the time needed to recover the file
$i$ knowing $x_i$. By stationarity, we can take $x_i=0$. Thus we
consider a file of size $l$
 which we store at time $t$ at location $0$ on the hardware whose free space space is equal to $\R(t)$.  The
first point  (resp. the last point) of the hardware occupied for the storage of this file
is equal to $d(t)$ (resp. to $d(t)+\taup^{(t)}_l$). This gives the distribution of the extremities of the portion of
the hardware used for the storage of a file.
$\newline$

Before the proofs, we make  some useful observations and give examples. First, we have  for every $\rho\geq 0$  (use (\ref{kappa})),
\be
\label{kappat}
\kappa ^{(t)}(\rho)= \rho  +\int_0^{\infty} (1-e^{-\rho x}) \Pi^{(t)} (\t{d}x)
\ee
and using (\ref{intPi})
\be
\label{relmes}
\bar{\Pi}^{(t)}(0)=t\bar{\nu}(0), \qquad  \int_0^{\infty} x \Pi^{(t)}(\t{d}x)  =\frac{\t{m}t}{1-\t{m}t}.
\ee
Using (\ref{egalintro}), we have also the following identity  of measures  on $[0,\infty[\times [0,\infty[$
\be
\label{egal}
x\P(\taup^{(t)}_l \in \t{d}x)\t{d}l=l\P(-Y^{(t)}_x \in \t{d}l)\t{d}x.
\ee
Finally, we give the distribution of the extremities of $\B$ :
\be
\label{loigd}
\P(-g(t) \in \t{d}x)=\P(d(t) \in \t{d}x)=(1-\t{m}t)\big(\delta _0 (\t{d}x)+\ind _{\{x>0\}} \bar{\Pi}^{(t)}(x) \t{d}x \big).
\ee
$\newline$

Let us consider three explicit examples
\begin{ex}
$(1) $ The basic example is when   $\nu=\delta _1$
(all files have  the same unit size as in the original parking problem in $\cite{chl}$). Then for all $x\in \RRR_+$ and $n
\in \N$,
\bea
&&\P(Y_x^{(t)}+x=n)=e^{-tx}\frac{(tx)^n}{n!}, \nonumber \\
\label{dirac}
&&\P(\taup^{(t)}_x=x+n)=\frac{x}{x+n}e^{-t(x+n)}\frac{(t(n+x))^{n}}{n!},
\eea
where the second identity follows from integrating $(\ref{egal})$ on $\{(x,l) : l \in [z,z+h], x-z=n\}$ and letting $h$ tend to $0$. Then,
$$\Pi^{(t)}(n)= \frac{(tn)^n}{n .n!} e^{-tn}$$ and $l(t)$ follows a
size biased  Borel law :
$$\P(l(t)=n)=(1-t) \frac{(tn)^n}{n!} e^{-tn}.$$

$(2) $ An other example where calculus can be made explicitly is the gamma case when $\nu(\t{d}l)=\ind _{\{l\geq 0\}}l^{-1}e^{-l}\t{d}l$. Note
that $\bar{\nu}(0)=\infty$ and
$\t{m}=1$. Then, for every $x\in\RRR_+$,

$$
\P(Y^{(t)}_x \in \t{d}z)=\ind _{[-x,\infty[}(z)\Gamma(tx)^{-1}e^{-(z+x)}(z+x)^{tx-1}\t{d}z,
$$
\be
\label{gamma}
\P(\taup^{(t)}_x\in \t{d}z)=\ind_{[x,\infty[}(z)x(z\Gamma(tz))^{-1}e^{-(z-x)}(z-x)^{tz-1}\t{d}z.
\ee
Further
$$\Pi^{(t)}(\t{d}z)= (z\Gamma(tz))^{-1}e^{-z}z^{tz-1}\t{d}z$$
and
$$\P(l(t) \in \t{d}x)=(1-t)\big(\delta_0(\t{d}x)+\Gamma(tz)^{-1}e^{-x}x^{tx-1}\t{d}x).$$
$\newline$

$(3) $ For the exponential distribution $\nu(\t{d}l)=\ind _{\{l\geq 0\}}e^{-l}\t{d}l$, we can get :

$$\Psi^{(t)}(\lambda)=\lambda (-1+\frac{t}{\lambda +1}), \quad \kappa^{(t)}(\lambda)=(\lambda+t-1+\sqrt{(\lambda+t-1)^2+4\lambda})/2.$$
\end{ex}
$\newline$

Finally, we specify two distributions involved in the storage  of the data. \\ \\
Writing $-g(t)=\gamma(t)$ (see (\ref{gam}) and (\ref{gsup})) and
using the identity of fluctuation ($\ref{fluc}$) gives an other expression for the Laplace transform
of $g(t)$ : For all $t\in [0,\t{1/m}[$ and  $\lambda \geq 0$, we have
\be
\label{gg}
\E\big(\t{exp}\big(\lambda g(t) \big) \big)=\t{exp}\bigg(\int_0^{\infty} (e^{-\lambda x}-1) x^{-1} \P(Y^{(t)}_x> 0)\t{d}x\bigg).
\ee
As a consequence, we see that the law of $g(t)$ is infinitively divisible. Moreover this expression will be useful to study the process $(g(t))_{t\in[0,1/\t{m}[}$
in \cite{vbb}. \\ \\
The quantity of data over $0$,  $R^{(t)}_0$  (see Section 2.1), is an increasing process equal to $(-I^{(t)}_0)_{t\geq0}$.
Its law is given by $S(t)$ (see (\ref{gam})) and, by $(\ref{sup2})$, its Laplace transform is then equal to
$$\lambda \longrightarrow \frac{(1-\t{m}t)\lambda }{\Psi^{(t)}(\lambda)}.$$
$\newline$

\subsection{Proofs}

\begin{proof}[Proof of Proposition \ref{dico}]
First $\t{m} < \infty$  entails that $\forall L\geq 0, \ \sum_{t_i\leq t, x_i \in [-L,L]} l_i< \infty$ and
condition $(\ref{condHH})$ is satisfied a.s. \\ \\
\ \ $\bullet$ \ If $t<\t{1/m}$, then $\E(Y^{(t)}_{-1})=1-\t{m}t>0$
and  the càdlàg version of
$(Y^{(t)}_{(-x)^-})_{x\geq 0}$ 
 is a Lévy process. So we have (see \cite{lev} Corollary 2 on page 190) :
 $$Y^{(t)}_x \stackrel{ \scriptscriptstyle x \rightarrow -\infty } {\longrightarrow} \infty \ \ \ \t{a.s.}$$
Then Proposition $\ref{disc}$ ensures that for every $t<\t{1/m}$,   $\C (t)= \{ x\in \RRR : \ Y^{(t)}_x > I^{(t)}_x\}\ne \RRR$ a.s. \\ \\
\ \ $\bullet$ \ If $t\geq \t{1/m}$, then  $\E(Y^{(t)}_{-1})\leq 0$ ensures (see \cite{lev} Corollary 2 on page 190) :
$$Y^{(t)}_x \stackrel{ \scriptscriptstyle x \rightarrow -\infty } {\longrightarrow} -\infty \  \t{a.s}\ \ \ \t{or} \ \ \ (Y^{(t)}_x)_{x\leq 0} \ \t{oscillates a.s in} -\infty $$
Similarly, we get that for every $t\geq \t{1/m}$, $\C (t)=\RRR$ a.s.
\end{proof}
$\newline$

For the other proofs, we fix  $t \in [0,\t{1/m}[$, which is omitted from the
notation of  processes for simplicity. \\

To prove the next proposition and the theorem, we need to establish first a regeneration property at the right extremities of the data blocks. In
that view, we
consider for every $x\geq 0$, the files arrived at the left/at the right of $x$ before time $t$ :
$$ \mathcal{P}_x:=\{(t_i,x_i,l_i) : t_i\leq t, \ x_i\leq x\}, \quad \mathcal{P}^x:=\{(t_i,x_i-x,l_i) : t_i\leq t, \ x_i> x\}.$$
\begin{Lem}
\label{reg} For every $x\geq 0$, $\mathcal{P}^{d_x(\R(t))}$ is
independent of $\mathcal{P}_{d_x(\R(t))}$ and distributed as
$\mathcal{P}^{0} $.
\end{Lem}
\begin{proof}
The simple Markov property for PPP states that for every $x\in\RRR$, $\mathcal{P}^{x}$ is
independent  of $\mathcal{P}_{x}$ and distributed as
$\mathcal{P}^{0}$. Clearly this extends to simple stopping times in the filtration $\sigma\big(\mathcal{P}_x\big)_{x\in\RRR}$
 and further to any stopping time in this filtration using the classical argument of approximation of stopping times
by a decreasing sequence of simple stopping times (see also \cite{mey}). As  $d_x(\R(t))$ is a
 stopping time in this filtration,  $\mathcal{P}^{d_x(\R(t))}$ is independent of
$\mathcal{P}_{d_x(\R(t))}$ and distributed as $\mathcal{P}^{0}$.
\end{proof}
$\\$
\begin{proof}[Proof of Proposition \ref{prob}]
$\bullet$ \ The free space at the right of $d_x(\R(t))$ at time $t$ is given by the point process of files arrived at the right of $d_x(\R(t))$ before time $t$. That
is, there exists a measurable functional $F$  such that for all $x \in \RRR$,
$$(\R(t)-d_x(\R(t)))\cap [0,\infty[=F\big(\mathcal{P}^{d_x(\R(t))}\big).$$
Similarly $\R(t)\cap ]-\infty,x]$ is $\mathcal{P}_{d_x(\R(t))}$
measurable. The previous lemma ensures then that
$(\R(t)-d_x(\R(t)))\cap [0,\infty[$ is independent of $\R(t)\cap
]-\infty,x]$ and
is distributed as  $(\R-d_0(\R(t)))\cap [0,\infty[$. \\ \\
$\bullet$ \ The stationarity of $\C(t)$ should be plain from the construction
of the covering and the fact that the law of a PPP with intensity $\t{d}x\otimes\nu$ is invariant by translation of
the first coordinate.  Stationarity can  also be viewed as a consequence of regeneration and  $\t{inf} \ \R(t)=-\infty$ (see Remark (4.11) in $\cite{mais}$). \\ \\
$\bullet$ \ The symmetry of $\R(t)^{cl}$ is  a consequence of the regeneration property and stationarity (see Lemma 6.5 in \cite{Tak} or
Corollary (7.19) in $\cite{Takk}$). \\ \\
$\bullet$  \ As a consequence of stationarity, $\P(x\in \C (t))$
does not depend on $x$ and  is equal to $\P(0\in \C (t))$. Following
Section 2.1, we write  $R_x:=Y_x-I_x$ the quantity of data over
$x$ so that the quantity of data
 stored in $[-L,L]$ is given for every $L>0$ by
$$ \mid \C (t) \cap [-L,L] \mid =R_{-L}+\big( \sum_{t_i\leq t, \ x_i \in ]-L,L]} l_i \big)-R_{L}.$$
By invariance of the PPP $\{(t_i,x_i,l_i) : i \in \N\}$ by translation of the second coordinate,
$$\P \big((2L)^{-1}R_{L} \geq \epsilon \big)=\P \big((2L)^{-1}R_{-L}\geq \epsilon \big)=\P\big ((2L)^{-1}R_{0}  \geq \epsilon \big) \stackrel{\scriptscriptstyle L \rightarrow \infty} {\longrightarrow} 0.$$
$$\t{Moreover using} \  (\ref{Lap1}), \    (2L)^{-1}\sum_{t_i\leq t, \ x_i \in ]-L,L[} l_i  \stackrel{ \scriptscriptstyle L \rightarrow \infty } {\longrightarrow} \t{m}t \quad \t{in probability}. \ \t{So}  \qquad \qquad \qquad \qquad \qquad $$
$$\E \big( (2L)^{-1} \mid \C (t) \cap [-L,L] \mid\big)\stackrel{ \scriptscriptstyle L \rightarrow \infty } {\longrightarrow} \t{m}t $$
and we conclude  with
$$ \E \big( \mid \C (t) \cap [-L,L] \mid\big)=\E\big(  \int_{-L}^{L} \ind _{\{x \in \C (t)\}}\t{d}x\big)=\int_{-L}^{L} \P(x\in \C (t))\t{d}x=2L\P(0\in\C(t)).$$
One can also give  a formal argument using Theorem  1 in $\cite{Tak}$ or $\P(0\in\C(t))=\P(l(t)>0)$ and Theorem $\ref{loi}$.
\end{proof}

$\newline$

\begin{proof}[Proof of  Theorem \ref{loi}]

(i) By symmetry of $\R(t)^{cl}$,  $\Rpt$ and  $\Rmt$ are  identically distributed. The
 regeneration property ensures that $\Rpt$ is independent of $(\Rmt,g(t),d(t))$. By
symmetry, $\Rmt$  is independent of $(g(t),d(t),\Rpt)$. So
$\Rpt$, $\Rmt$ and $(g(t),d(t))$
are independent.  \\ \\

(ii) As $\Rp(t)$ is a.s. the union of intervals of the form $[a,b[$, then for
every $x\geq 0$,
$$\taup_{\vert\R(t)\cap [0,x]\vert}=d_x(\R(t)), \qquad \taup_x=d_{\taup_x}(\R(t)) \qquad  a.s.$$
So  the range of $\taup$ is equal to $\Rpt$. The fact $\taup$ is a subordinator will be proved below but could be derived now
from the  regeneration property of $\Rp(t)$. Similarly the range of
$\taum$ is equal to $\Rmt$. \\
Moreover,  $\t{d}Y=-1$ on $\R(t)$ and $Y_{a-}=Y_b$ if $[a,b[$ is an
interval component of $\C(t)$. By integrating on $[d(t),d(t)+y]$,
we have a.s  for every $y\geq 0$ such that $d(t)+y\in\R(t)$,
$$Y_{y+d(t)}-Y_y=- \vert\R(t)\cap [d(t),d(t)+y]\vert.$$
Then using again the definition of $\taup$ given in Section 3.1 and that $\Rpt$ is the range of $\taup$,
\bea
\taup_x&=& \t{inf}\{y\geq 0 : y \in \Rpt, \  \vert\Rpt\cap [0,y]\vert>x\} \nonumber \\
&=&\t{inf}\{y\geq 0 : d(t)+ y \in \R(t), \ \vert  \R(t)\cap [d(t),d(t)+y]\vert> x\} \nonumber \\
\label{taup}
&=&\t{inf}\{y\geq 0 : Y_{y+d(t)}-Y_{d(t)}<-x\}
\eea
Moreover,
$$Y_{y+d(t)}-Y_{d(t)}=-y+\sum_{\substack{(t_i,x_i,l_i)\in\mathcal{P}^{d(t)}\\ 0\leq x_i\leq y}} l_i$$
and Lemma \ref{reg} entails that $\mathcal{P}^{d(t)}$ is distributed as a PPP on $[0,t]\times \RRR_+\times \RRR_+$ with
intensity
$\t{d}s\otimes \t{d}x\otimes \nu(\t{d}l)$. So
$(Y_{y+d(t)}-Y_{d(t)})_{y\geq 0}$ is a Lévy process with bounded variation and drift $-1$ which verifies condition  (\ref{condHH})
(use (\ref{lim2}) and $-1+\t{m}t<0$). Then  Theorem $\ref{tau}$ entails that $\taup$ is a subordinator
whose  Laplace exponent  is the inverse function of $-\Psi^{(t)}$. \\ \\

As $\Rmt$ is distributed as $\Rpt$,  $\taum$ is distributed as  $\taup$ by definition. \\ \\

(iii) We determine now the distribution of $(g(t),d(t))$ using
fluctuation theory, which enables us  to get identities useful for the
rest of the work. We write $(\w{Y}_x)_{x\geq 0}$ for the càdlàg
version of $(-Y_{-x})_{x\geq 0}$ and \be \label{gam}
S{(t)}:=\t{sup}\{\w{Y}_x,x\geq 0\}=-I_0, \qquad
\gamma{(t)}:=\t{arg}(S(t))=\inf\{x\geq 0 : \w{Y}_{x}=S(t)\}. \ee
 Using (\ref{comp}) and the fact
that $Y$ has no negative jumps, we have \bea
 g(t)&=&g_0(\R(t))= \t{sup}\{x\leq 0 : \ Y_x=I_x\} \nonumber \\
&=&\t{sup}\{x\leq 0 : \ Y_{x-}=I_0\}=-\t{inf}\{x\geq 0 : \ \w{Y}_x=-I_0\}\nonumber \\
\label{gsup}
&=&-\gamma(t)
\eea
Using again (\ref{comp}) and the fact that $(Y_x)_{x\geq 0}$ is regular for $]-\infty,0[$ (see  \cite{lev} Proposition 8
on page 84), we have also a.s.
\bea
d(t)&=&\t{inf}\{x > 0 :  Y_x=I_x\}=\t{inf}\{x> 0 : Y_x=I_0\} \nonumber  \\
&=&\t{inf}\{x> 0 : Y_x<I_0\}=\t{inf}\{x> 0 : Y_x<-S(t)\} =T_{S(t)}  \nonumber
\eea
where $(T_x)_{x \geq 0}$ is distributed as $(\taup_x)_{x\geq 0}$ by (\ref{taup}) and $(T_x)_{x \geq 0}$ is independent of
$(S(t),\gamma(t))$ since $(Y_x)_{x\geq 0}$ is independent of $(Y_x)_{x\leq 0}$.
 Then for all  $\lambda ,\mu \geq  0$ with $\lambda\ne\mu$ :
\bea
\E\big(\t{exp}(\lambda g(t)  -\mu d(t))\big)&=& \E\big(\t{exp}(-\lambda \gamma(t))\E(\t{exp}(-\mu T_{S(t)}))\big) \nonumber \\
&=&\E\big(\t{exp}(-\lambda \gamma(t)-\kappa^{(t)}(\mu)S(t))\big) \nonumber \\
\label{lapdoubb}
&=& -[\Psi^{(t)}]'(0) \frac{\kappa^{(t)}(\lambda)-\kappa^{(t)}(\mu)}{\lambda-\mu}  \ \ \ \ \t{using} \ (\ref{sup}) \\
\label{lapdoub}
&=&(1-\t{m}t) \frac{\kappa^{(t)}(\lambda)-\kappa^{(t)}(\mu)}{\lambda-\mu} \ \ \ \ \  \ \t{using} \ (\ref{lim2})
\eea
which gives the distributions of $d(t)$, $g(t)$ and $l(t)$ letting respectively $\lambda=0$, $\mu=0$ and $\lambda\rightarrow \mu$. Computing
then the Laplace transform  of $(-Ul(t),(1-U)l(t))$ where $U$ is a uniform
random variable on $[0,1]$ independent of $l(t)$ gives the right hand side of ($\ref{lapdoub}$).  So
 $(g(t),d(t))=(-U'l(t),(1-U')l(t))$, where
 $U'$ is a uniform
random variable on $[0,1]$ independent of $l(t)$ .
\end{proof}

\begin{Rque}
We have proved above that $\Rmt$ is distributed as $\Rpt$, which entails that the  last passage-time-process of
 the post-infimum  process of $(-Y_x)_{x\geq 0}$ is distributed as
the first-passage-time process of   $(-Y_x)_{x\geq 0}$. \\
This result is also a consequence of
 the fact that  the post-infimum process of $(-Y_x)_{x\geq 0}$ is distributed
as the Lévy process $(-Y_x)_{x\geq 0}$ conditioned to stay
positive \cite{Mill}, whose last-passage-time process is a
subordinator with Laplace exponent  $\kappa$ (see  Exercise 3 on
page 213 in \cite{lev}).
\end{Rque}
$\newline$

\begin{proof}[Proof of  Corollary \ref{rq}]
As $\bar{\nu}(0)<\infty$, then $\bar{\Pi}(0)=t\bar{\nu}(0)<\infty$ (see (\ref{relmes})). So $\taup$ is the sum of a drift
and a compound Poisson process. That is, there exists a Poisson process $(N_x)_{x\geq 0}$ of intensity $t\bar{\nu}(0)$
 and a sequence $(X_i)_{i\in\N}$ of iid variables of law $\nu/\bar{\nu}(0)$ independent of $(N_x)_{x\geq 0}$ such that
$$\taup_x=x+\sum_{i=1}^{N_x} X_i, \qquad x\geq 0.$$
As $\Rpt$ is the range of $\taup$, the number of data blocks of $\C(t)$ between $d(t)$ and $d(t)+\taup_x$ is  equal to
the number of jumps of $\taup$ before $x$, that is $N_x$. Thus,
$$\frac{\t{number of data blocks in} \ [d(t),d(t)+\taup_x]}{\taup_x}= \frac{N_x}{\taup_x}\stackrel{x\rightarrow \infty}{\longrightarrow } \frac{\E(N_1)}{\E(\taup_1)}=t\bar{\nu}(0)(1-\t{m}t) \ \ \t{a.s.}$$
by the law of large numbers (see \cite{lev} on page 92). This completes the proof.
\end{proof}
$\newline$

\section{Asymptotics  at saturation of the hardware}

We focus now on the asymptotic behavior of $\R(t)$  when $t$ tends
to $\t{1/m}$, that is when the hardware is becoming full.  First, note that if $\nu$ has a finite second moment,  then
$$\E \big( l(t) \big)=\frac{\int_0^{\infty} l^2
\nu(\t{d}l)}{(1-\t{m}t)^2}.$$
Thus we may expect  that if $\nu$ has a finite second moment, then
$(1-\t{m}t)^{2}l(t)$ should converge in distribution as $t$ tends
to $1/\t{m}$. Indeed, in the particular case $\nu=\delta_1$ or in
the conditions of Corollary 2.4 in \cite{mono}, we have an
expression of $\Pi^{(t)}(\t{d}x)$ and we can prove that
$(1-\t{m}t)^{2}l(t)$ does converge in distribution to a gamma
variable. More generally, we shall prove that  the rescaled
free space $(1-\t{m}t)^{2}\R(t)$ converges in distribution  as
$t$ tends to $1/\t{m}$. In that view, we need to prove  that the
process $(Y^{(t)}_{(1-\t{m}t)^{-2}x})_{x\in \RRR}$ converges after
suitable rescaling to a random process. Thanks to $(\ref{comp})$,
$(1-\t{m}t)^2\R(t)$ should then converge to the set of points where
this limiting process coincides with its infimum process. We shall
also han\t{d}le the case where $\nu$ has an infinite second moment and
find the correct normalization.
 \\
 
  In queuing systems, asymptotics at saturation are known as heavy traffic approximation ($\rho=t\t{m}\rightarrow 1$), which
depend similarly on the tail of $\nu$. And for  $\nu$ finite, results given here could be directly
derived from results in queuing theory  (See III.7.2 in
 \cite{cohen} or \cite{King} if $\nu$ has a second moment order and \cite{boxma}  for heavy tail of $\nu$). The main
difference here is that $\nu$ can be infinite and we  consider the whole random set of occupied
locations. Moreover,  as explained below, asymptotics of $\R(t)$  can not be directly derived from asymptotics of  $Y$ or the workload $R$. \\

Following the notation in $\cite{bem}$, we say that $\nu \in \D_{2+}$ if $\nu$ has a finite second moment 
$\t{m}_2:=\int_0^{\infty}l^2\nu(\t{d}l)$. For $\alpha \in ]1,2]$, we say that $\nu \in \D_{\alpha}$ whenever
$$\exists C>0 \ \ \t{such that} \ \ \bar{\nu}(x)\stackrel{x\rightarrow \infty}{\sim} Cx^{-\alpha}$$
and we put for $\alpha \in ]1,2[$ :
$$C_{\alpha}:=\bigg(\frac{C\Gamma(2-\alpha)}{\t{m}_2(\alpha-1)}\bigg)^{1/\alpha}.$$
We denote by $(B_z)_{z\in\RRR}$ a two-sided Brownian motion,
i.e.  $(B_x)_{x\geq 0}$ and $(B_{-x})_{x\geq 0}$  are independent
standard Brownian motions. For $\alpha \in ]1,2[$ , we denote by
$(\sigma^{(\alpha)}_z)_{z\in\RRR}$ a càdlàg process  with independent and stationary
increments such that  $(\sigma^{(\alpha)}_x)_{x\geq 0}$ is  a
standard spectrally positive stable Lévy process with index
$\alpha$ :
$$\forall \ x\geq 0, \lambda \geq 0, \ \ \E(\t{exp}(-\lambda \sigma^{(\alpha)}_x)=\t{exp}(x\lambda^{\alpha}).$$
$\newline$

We introduce now  the following functions and processes defined for all $(t,x,z) \in [0,1/\t{m}[\times \RRR^*_+\times\RRR$ and $\alpha \in ]1,2[$ by \\
\bea
\epsilon_{2+}(t)=(1-\t{m}t)^2 \qquad \  \ && f_{2+}(x)=1/\sqrt{x}  \qquad   \ \ \ \ \  Y^{2+,\lambda}_z=-\lambda z+ \sqrt{\t{m}_2/\t{m}} B_z \nonumber \\
\epsilon_{2}(t)=2\frac{(1-\t{m}t)^{2}}{-\t{log}((1-\t{m}t))} \  &&  f_2(x)=\sqrt{\t{log}(x)/x} \qquad  Y^{2,\lambda}_z=-\lambda z+\sqrt{C/\t{m}} B_z \nonumber \\
\epsilon_{\alpha}(t)=(1-\t{m}t)^{\frac{\alpha}{\alpha-1}} \ \qquad && f_{\alpha}(x)=x^{1/\alpha-1} \qquad \ \ \ \ \ Y^{\alpha,\lambda}_z=-\lambda z+C_{\alpha} \sigma^{(\alpha)}_z \nonumber
\eea
and the infimum process defined for $x\in\RRR$ by $I^{\alpha,\lambda}_x:=\t{inf}\{Y^{\alpha,\lambda}_y : y \leq x\}$. \\ \\

We have the following weak convergence result for the Hausdorff
metric defined in Section 2.

\begin{Thm}
\label{asympt}
If $\nu \in \D_{\alpha}$ ($\alpha \in [1,2]\cup \{2+\}$), then $\epsilon_{\alpha}(t). \R(t)^{cl}$ converges
weakly in $\H(\RRR)$ as $t$
tends to $1/\emph{m}$ to $\{x\in \RRR : Y^{\alpha,1}_x=I^{\alpha,1}_x\}^{cl}$ .
\end{Thm}
First we prove the convergence of the Laplace exponent $\Psi^{(t)}$ after suitable rescaling  as $t$ tends to $1/\t{m}$, which
ensures the convergence of   the Lévy process $Y^{(t)}$ after suitable rescaling  (see Lemma \ref{process}).
These convergences  will not a priori entail the convergence of the random set $\epsilon_{\alpha}(t). \R^{cl}(t)$
since they do not entail
the convergence of excursions. Nevertheless, they will entail the convergence of $\kappa^{(t)}$ since
$\kappa^{(t)}\circ (-\Psi^{(t)})=\rm{Id}$ (Lemma \ref{processtau}).
Then we get the convergence of $\tau^{(t)}$   as $t$ tends to infinity and thus of its range
$\epsilon_{\alpha}(t). \R^{cl}(t)$.
$\newline$
\begin{Rque}
More generally, if  $\bar{\nu}$ is regularly varying   at infinity  with  index $-\alpha \in ]-1,-2[$, then
we have the following weak convergence in $\H(\RRR)$
$$x^{-1}\R ((1-x\bar{\nu}(x))/\t{m})^{cl} \stackrel{x\rightarrow \infty}{\Longrightarrow }
\{x\in \RRR : Y^{\alpha,1}_x=I^{\alpha,1}_x\}^{cl}\quad  \t{with} \ C=1.$$
For instance, the case $\bar{\nu}(x) \stackrel{x\rightarrow
\infty}{\sim} cx^{-\alpha} \t{log}(x)^{\beta}$ with
$(\alpha,\beta,c) \in ]1,2[\times\RRR\times  \RRR^*_+$ leads to
$$\big((1-\t{m}t)\t{log}(1/(1-\t{m}t))^{-\beta}\big)^{\frac{1}{\alpha-1}} \R(t)^{cl}\stackrel{t\rightarrow 1/\t{m}}{\Longrightarrow }
\{x\in \RRR : Y^{\alpha,1}_x=I^{\alpha,1}_x\}^{cl}   \ \ \  \t{with} \ C=c/(\alpha-1)^{\beta}.$$
If $\bar{\nu}$  is regularly varying at infinity with index  $-2$, there are many cases to consider.  \\
\end{Rque}
$\newline$
We get then the asymptotic of $(g(t),d(t))$ :
\begin{Cor}
\label{asymptl}
If $\nu \in \D_{\alpha}$ ($\alpha \in [1,2]\cup \{2+\}$), then
$\epsilon_{\alpha}(t).(g(t),d(t))$ converges weakly as $t$ tends to $1/\emph{m}$ to
$(\emph{sup}\{x\leq 0 : Y^{\alpha,1}_x=I^{\alpha,1}_0\}, \emph{inf}\{x\geq 0 : Y^{\alpha,1}_x=I^{\alpha,1}_0\})$. \\
If $\nu \in \D_{2+}$ (resp. $\D_{2}$), $\epsilon_{\alpha}(t). l(t)$ converges weakly to  a gamma variable   with parameter $(1/2,\emph{m}/(4\emph{m}_2))$ (resp.
$(1/2,\emph{m}/4)$). \\
\end{Cor}
\begin{Rque}  The density of data blocks of size $\t{d}x$ in  $\epsilon_{\alpha}(t). \R(t)^{cl}$
 is equal to $\frac{\t{m}t}{1-\t{m}t} \Pi^{(t)}(\t{d}x)$. By the previous theorem or corollary, this density
 converges weakly as $t$ tends to $1/\t{m}$ to the density of data block of size $\t{d}x$ of the limit
covering  $\{x\in \RRR : Y^{\alpha,1}_x=I^{\alpha,1}_x\}^{cl}$.
This limit density, denoted by $\Pi^{\alpha,1}(\t{d}x)$,
  can be computed explicitely in the cases  $\nu \in \D_{\alpha} \ (\alpha \in \{2,2+\})$, thanks to the last corollary :
$$\Pi^{2+,1}(\t{d}x)=\sqrt{\frac{\t{m}}{4\pi\t{m}_2 x^3}} \t{exp}\big(-\frac{\t{m}}{4\t{m}_2}x\big), \qquad \Pi^{2,1}(\t{d}x)=\sqrt{\frac{\t{m}}{4\pi x^3}} \t{exp}\big(-\frac{\t{m}}{4}x\big).$$
Note that is also the Lévy measure of the limit covering  $\{x\in
\RRR : Y^{\alpha,1}_x=I^{\alpha,1}_x\}^{cl}$.
\end{Rque}
$\newline$
$\newline$

If we look at $\C(t)$ in a window of size $x$ and let $x$ tend to infinity, we observe :
\begin{Thm}
\label{thtrst}
If $\nu \in \D_{\alpha}$ ($\alpha \in [1,2]\cup \{2+\}$),  $x$ tends to infinity and $t$ to $1/\emph{m}$
such that $1-\emph{m}t\sim \lambda f_{\alpha}(x)$ with $\lambda>0$, then
$x^{-1} (\R(t)^{cl}\cap [0,x])$ converges weakly in $\H([0,1])$ to $\{x \in [0,1] : Y^{\alpha,\lambda}_x=I^{\alpha,\lambda}_x\}^{cl}$.
\end{Thm}
$\newline$
Thus as in $\cite{chl}$, we observe a  phase transition of   the size of largest block of data in $[0,x]$ as $x\rightarrow \infty$
 according
to the rate of filling of the hardware. More precisely, denoting \\
   $B_1(x,t)=\mid I_1(x,t)\mid$ where $(I_j(x,t))_{j\geq 1}$ is the sequence of component intervals
of $\C(t)\cap [0,x]$ ranked by decreasing order of size, we have :

\begin{Cor}
\label{cortrst}
Let $\nu \in \D_{\alpha}$ ($\alpha \in [1,2]\cup \{2+\}$),    $x$ tend to infinity and $t$ to $1/\emph{m}$ : \\
- If \ $1-\emph{m}t\sim \lambda f_{\alpha}(x)$ with $\lambda>0$, then $B_1(x,t)/x$ converges in distribution to  the largest length of excursion
of $(Y^{\alpha,\lambda}_x-I^{\alpha,\lambda}_x)_{x\in [0,1]}$.\\
- If \ $f_{\alpha}(x)=o(1-\emph{m}t)$, then $B_1(x,t)/x\stackrel{\P}{\longrightarrow }0$ . \\
- If \ $1-\emph{m}t=o(f_{\alpha}(x))$, then $B_1(x,t)/x\stackrel{\P}{\longrightarrow }1$.
\end{Cor}

The phase transition occurs at time $t$ such that $1-\t{m}t\sim \lambda f_{\alpha}(x)$ with $\lambda>0$. The more data
arrive in small files (i.e. the faster $\bar{\nu}(x)$  tends  to zero as $x$ tends to infinity), the later
 the phase transition occurs.
The phase transition in $\cite{chl}$ or $\cite{bem}$ uses the bridges of the processes involved here. A
consequence is that  in our model, $B_1(t,x)/x$  tends to zero or one with a positive probability at phase transition, which
is not the case for the parking problem in  $\cite{chl}$ or $\cite{bem}$. More precisely,
denoting by $B_{\alpha,\lambda}$
 the law of the largest length of excursion
of $(Y^{\alpha,\lambda}_x-I^{\alpha,\lambda}_x)_{x\in [0,1]}$, we have :
$$\forall \ (\lambda,\alpha)\in \RRR_+^*\times ]1,2[\cup\{2+\}, \qquad\P(B_{\alpha,\lambda}=0)>0, \qquad \P(B_{\alpha,\lambda}=1)>0.$$

$\newline$

For the proofs of the theorems,  we introduce $\Psi^{\alpha,\lambda}$ the Laplace exponent (see (\ref{lapdef})) of $Y^{\alpha,\lambda}$
given  for $y\geq 0$, $\lambda\geq 0$ and $\alpha \in ]1,2[$  by
\be
\label{lapdrift}
\Psi^{2+,\lambda}(y)= -\lambda y -\frac{\t{m}_2}{\t{m}} \frac{y^2}{2}, \qquad
\Psi^{2,\lambda}(y)=  -\lambda y -\frac{C}{\t{m}} \frac{y^2}{2}, \qquad
\Psi^{\alpha,\lambda}(y)=  -\lambda y - (C_{\alpha}y)^{\alpha}. \nonumber \\
\ee We denote by $\DD$  the space of càdlàg function from $\RRR_+$
to $\RRR$ which we
endow with the Skorokhod topology (see \cite{jac} on page 292). First, we prove the weak convergence of   $Y^{(t)}$ after suitable rescaling. \\

\begin{Lem}\label{process}  If  $\nu \in \D_{\alpha}$ ($\alpha \in [1,2]\cup \{2+\}$), then for all  $y\geq 0$ and $\lambda> 0$ :  \\
$$\epsilon_{\alpha}(t)^{-1} \Psi^{(t)}(\epsilon_{\alpha}(t)(1-\emph{m}t)^{-1}y) \stackrel{t\rightarrow 1/\emph{m}}{\longrightarrow }  \Psi^{\alpha,1}(y),$$
$$x\Psi^{((1-\lambda
f_{\alpha}(x))/\emph{m})}((xf_{\alpha}(x))^{-1}y)\stackrel{x\rightarrow
\infty}{\longrightarrow }  \Psi^{\alpha,\lambda}(y),$$
which entail the following weak convergences of processes in $\DD$ :
$$ \big( \epsilon_{\alpha}(t)(1-\emph{m}t)^{-1} Y_{\epsilon_{\alpha}(t)^{-1}y}^{(t)}\big)_{y\geq 0} \stackrel{t\rightarrow 1/\emph{m}}{\Longrightarrow  } (Y^{\alpha,1}_y)_{y\geq 0},
$$
$$\big((xf_{\alpha}(x))^{-1} Y_{xy}^{((1-\lambda
f_{\alpha}(x))/\emph{m})}\big)_{y\geq 0}\stackrel{x\rightarrow \infty}{\Longrightarrow } (Y^{\alpha,\lambda}_y)_{y\geq 0}.$$
\end{Lem}

\begin{Rque}
If  $\bar{\nu}$ is regularly varying at infinity with index $-\alpha \in ]-1,-2[$,  then
$\bar{\nu}(x)^{-1}\Psi^{((1-\lambda x\bar{\nu}(x))/\t{m})} (x^{-1}y)$
converges  to $\Psi^{\alpha,\lambda}(y)$ as $x$ tends to infinity.  \\
\end{Rque}
\begin{proof}[Proof of Lemma \ref{process}] Using (\ref{Lap2}), we have
\be
\label{laptrs}
x\Psi^{(t)}(y)=xy(\t{m}t-1-t\int_0^{\infty}(1-e^{-y u})\bar{\nu}(u)\t{d}u).
\ee
We   han\t{d}le now the different cases : \\ \\
$\bullet$  Case  $\nu \in \D_{2+}$.   Using $\vert 1-e^{-yu}\vert/y\leq u \ \ (u\geq 0)$ and dominated
convergence theorem gives :
$$\int_0^{\infty}(1-e^{-yu})\bar{\nu}(u)\t{d}u \ \stackrel{y\rightarrow 0}{\sim} \ y \int_0^{\infty}u\bar{\nu}(u)\t{d}u=\frac{y\t{m}_2}{2}$$
which proves the first part of the lemma  using $(\ref{laptrs})$. \\ \\
$\bullet$ Case $\nu \in \D_{\alpha}$ with $\alpha \in ]1,2[$. Using that $(u/y)^{\alpha}\bar{\nu}(u/y)$ is
bounded, we apply dominated convergence theorem and  get
\bea
\label{lapch}
\int_0^{\infty}(1-e^{-y u})\bar{\nu}(u)\t{d}u &=&
 y^{-1}\int_0^{\infty}(1-e^{-u})\bar{\nu}(u/y)\t{d}u \\
&\stackrel{y\rightarrow 0}{\sim}& y^{-1} C(y^{-1})^{-\alpha}\int_0^{\infty}(1-e^{-u})u^{-\alpha}\t{d}u  \nonumber \\
&\stackrel{y\rightarrow 0}{\sim}& C \frac{\Gamma(2-\alpha)}{\alpha-1}y^{\alpha-1} \nonumber
\eea
which proves  the first part of the lemma using $(\ref{laptrs})$. \\ \\
$\bullet $ Case $\nu$ is regularly varying at infinity with index $-\alpha\in]-1,-2[$. First,
$$ \int_0^{1/\sqrt{y}} (1-e^{-yu})\bar{\nu}(u)\t{d}u\leq y \int_0^{1/\sqrt{y}}u\bar{\nu}(u)\t{d}u \stackrel{y\rightarrow 0}{\sim} y(1/\sqrt{y})^{2-\alpha}=y^{\alpha/2}$$
Moreover for every $u>0$, $\bar{\nu}(u/y)\stackrel{y\rightarrow 0}{\sim}\bar{\nu}(y) u^{-\alpha}$. Let $\delta>0$ such
 that  $-2<-\alpha-\delta<-\alpha+\delta<-1$. By Potter's theorem (page 25 in \cite{bin}) ensures
that for all $y$ small enough and  $u$ large enough,
$$\frac{\bar{\nu}(u/y)}{\bar{\nu}(1/y)}\leq 2\t{max}(u^{-\alpha+\delta},u^{-\alpha-\delta}).$$
So we can apply the  dominated convergence theorem to  get
$$ \int_0^{1/\sqrt{y}} (1-e^{-yu})\bar{\nu}(u)\t{d}u= y^{-1}\int_{\sqrt{y}}^{\infty}(1-e^{-u})\bar{\nu}(u/y)\t{d}u \stackrel{y\rightarrow 0}{\sim} \frac{\Gamma(2-\alpha)}{\alpha-1}  y^{-1}\bar{\nu}(1/y).$$
As $y^{\alpha/2}=o(y^{-1}\bar{\nu}(1/y))$ \ $(y\rightarrow 0)$, we can complete the proof with
$$\int_0^{\infty}(1-e^{-y u})\bar{\nu}(u)\t{d}u\stackrel{y\rightarrow 0}{\sim} \frac{\Gamma(2-\alpha)}{\alpha-1}  y^{-1}\bar{\nu}(1/y).$$
$\newline$
$\bullet$  Case  $\nu \in \D_{2}$. We split the integral. First, we have
\bea
\int_0^{1/\sqrt{y}}(1-e^{-y u})\bar{\nu}(u)\t{d}u  &\stackrel{y\rightarrow 0}{\sim}&  y\int_0^{1/\sqrt{y}} u\bar{\nu}(u)\t{d}u \ \  \stackrel{y\rightarrow 0}{\sim}   \ \ Cy \t{log}(1/y)/2. \nonumber
\eea
since $\int_0^{1/\sqrt{y}}(1-e^{-y u}+yu)\bar{\nu}(u)\t{d}u=o (y \int_0^{1/\sqrt{y}} u\bar{\nu}(u)\t{d}u)$. Moreover,
\bea
\int_{1/\sqrt{y}}^{\infty}(1-e^{-yu})\bar{\nu}(u)\t{d}u \ &=& \ y^{-1} \int_{\sqrt{y}}^{\infty}(1-e^{-u})\bar{\nu}(u/y)\t{d}u \nonumber \\
\ &\stackrel{y\rightarrow 0}{\sim}& \ Cy  \int_{\sqrt{y}}^{\infty}(1-e^{-u})u^{-2} \t{d}u \quad \t{using} \ \nu \in \D_{2} \nonumber \\
\ &\stackrel{y\rightarrow 0}{\sim}& \ Cy  \int_{\sqrt{y}}^{1} u^{-1} \t{d}u\ = \ Cy \t{log}(1/y)/2 \nonumber
\eea
Then
$$\int_0^{\infty}(1-e^{-y u})\bar{\nu}(u)\t{d}u \ \stackrel{y\rightarrow \infty}{\sim} \ Cy \t{log}(1/y)$$
which proves the first part of the lemma using $(\ref{laptrs})$.\\ \\

These  convergences ensure  the convergence of the
finite-dimensional distributions of the processes. The weak
convergence in $\DD$, which is the second part of the lemma,  follows
from Theorem 13.17 in \cite{kal}.
\end{proof}
$\newline$

In the spirit of  Section 3, we   introduce the expected limit set, that
is the free space of the covering associated with
$Y^{\alpha,\lambda}$,  and the extremities of the block containing
$0$.
$$ \R(\alpha,\lambda):=\{x \in \RRR : Y^{\alpha,\lambda}_x=I^{\alpha,\lambda}_x\},
\qquad g(\alpha,\lambda):=g_0( \R(\alpha,\lambda)), \qquad
d(\alpha,\lambda):=d_0(\R(\alpha,\lambda)).$$ We have the following analog
of  Theorem \ref{loi}. $\Rpa$ and   $\Rma$ are
independent,  identically distributed and independent of
$(g(\alpha,\lambda),d(\alpha,\lambda))$. Moreover  $\Rpa$ and
$\Rma$ are respectively the range  of the
  subordinators $\taup^{\alpha,\lambda}$ and $\taum^{\alpha,\lambda}$, whose  Laplace exponent $\kappa^{\alpha,\lambda}$  is the
inverse function of $-\Psi^{\alpha,\lambda}$. Finally, using $\big[\Psi^{\alpha,\lambda}\big]'(0)=-\lambda$, the
counterpart of $(\ref{lapdoubb})$ gives for  $\rho,\mu\geq 0$
and $\rho\ne\mu$  :
 \be \label{laplenght} \E\big(\t{exp}(\rho
g(\alpha,\lambda)  -\mu d(\alpha,\lambda))\big)=
\lambda\frac{\kappa^{\alpha,\lambda}(\rho)-\kappa^{\alpha,\lambda}(\mu)}{\rho-\mu}.
\ee \\
The proof of these results follow the proof of Section 3.2, except for two points : \\  \\
1) We cannot use  the point process of files to prove the stationarity and regeneration property of $\R(\alpha,\lambda)$ and
we must use the process $ Y^{\alpha,\lambda}$ instead. The stationarity is  a direct consequence of the stationarity of
$\big(Y^{\alpha,\lambda}_x-I^{\alpha,\lambda}_x\big)_{x\in\RRR}$.  The regeneration property is a consequence of the counterpart of Lemma \ref{reg}
which can be stated as follows. For all $x\in\RRR$,
$$\big(Y^{\alpha,\lambda}_{d_x( \R(\alpha,\lambda))+y}-Y^{\alpha,\lambda}_{d_x( \R(\alpha,\lambda))} \big)_{y\geq 0}  \ \
\t{is independent of}  \ \ \big(Y^{\alpha,\lambda}_{d_x( \R(\alpha,\lambda))-y}-Y^{\alpha,\lambda}_{d_x( \R(\alpha,\lambda))} \big)_{y\geq 0}$$
 and distributed as   $\big(Y^{\alpha,\lambda}_{y}\big)_{y\geq 0}$. As
Lemma \ref{reg}, this property is  an extension to  the stopping time $d_x( \R(\alpha,\lambda))$ of the following obvious result :
$\big(Y^{\alpha,\lambda}_{x+y}-Y^{\alpha,\lambda}_{x}\big)_{y\geq 0}$ is independent of
$\big(Y^{\alpha,\lambda}_{x-y}-Y^{\alpha,\lambda}_{x}\big)_{y\geq 0}$
and distributed as  $\big(Y^{\alpha,\lambda}_{y}\big)_{y\geq 0}$. \\ \\
2) It is convenient to define directly $(\taup^{\alpha,\lambda}_x)_{x\geq 0}$ by
$$\taup^{\alpha,\lambda}_x:=\t{inf}\{y\geq 0 : Y^{\alpha,\lambda}_{d(\alpha,\lambda)+y}-Y^{\alpha,\lambda}_{d(\alpha,\lambda)}<-x\}.$$
For  $\lambda>0$, $[\Psi^{\alpha,\lambda}] ' (0)=-\lambda<0$ so we can apply Theorem \ref{tau} and $\taup^{\alpha,\lambda}$
is a subordinator whose Laplace $\kappa^{\alpha,\lambda}$  is the
inverse function of $-\Psi^{\alpha,\lambda}$. Moreover its range is a.s. equal to $\Rpa$, since the Lévy process
$(Y^{\alpha,\lambda}_{d(\alpha,\lambda)+y}-Y^{\alpha,\lambda}_{d(\alpha,\lambda)})_{y\geq 0}$  is  regular for
$]-\infty,0[$ (Proposition 8 on page 84
in \cite{lev}).
$\newline$

To prove the theorems, we need a final lemma, which states the convergence of the Laplace exponent of $\Rpt$.

\begin{Lem}
\label{processtau}
If  $\nu \in \D_{\alpha}$ ($\alpha \in [1,2]\cup \{2+\}$), then for all  $z\geq 0$ and $\lambda> 0$,  \\
$$ (1-\emph{m}t)\epsilon_{\alpha}(t)^{-1}\kappa^{(t)}(\epsilon_{\alpha}(t) z) \stackrel{t\rightarrow 1/\emph{m}}{\longrightarrow } \kappa^{\alpha,1}(z),$$
$$xf_{\alpha}(x)\kappa^{((1-\lambda f_{\alpha}(x))/\emph{m})}(x^{-1}z) \stackrel{x\rightarrow \infty}{\longrightarrow }   \kappa^{\alpha,\lambda}(z).
$$
\end{Lem}
\begin{Rque}
If  $\bar{\nu}$ is regularly varying at infinity of index $-\alpha \in ]-1,-2[$, we have similarly
$$\bar{\nu}(x)^{-1}\kappa^{((1-\lambda x\bar{\nu}(x))/\t{m})}(x^{-1}z)\stackrel{x\rightarrow \infty}{\longrightarrow }  \kappa ^{\alpha,\lambda}(z).$$
\end{Rque}
\begin{proof}
First we prove that \be \label{equiv} \alpha(t) \
\stackrel{t\rightarrow 1/\t{m}}{\sim} \ \beta (t) \quad
\Rightarrow \quad \kappa^{(t)}(\alpha(t))  \
\stackrel{t\rightarrow 1/\t{m}}{\sim} \ \kappa^{(t)}(\beta (t)). \ee
Indeed the function $u \in \RRR^*_+\mapsto \frac{1-e^{-u}}{u}$
decreases so for all $x\geq 0$ and  $u,v>0$, we have :
$$ \t{min}(\frac{u}{v},1) \leq \frac{1-e^{-ux}}{1-e^{-vx}}\leq \t{max}(\frac{u}{v},1),$$
which gives
$$ \t{min}(\frac{\alpha(t)}{\beta(t)},1) \leq  \frac{\int_0^{\infty} (1-e^{-\alpha(t)x})\Pi^{(t)}(\t{d}x) }{ \int_0^{\infty} (1-e^{-\beta(t)x})\Pi^{(t)}(\t{d}x)} \leq  \t{max}(\frac{\alpha(t)}{\beta(t)},1)   $$
and proves  ($\ref{equiv}$) recalling $(\ref{kappat}$). \\

Then the first part of Lemma $\ref{process}$ and the identity $\kappa^{(t)}\circ (-\Psi^{(t)})=\rm{Id}$ give the first part of  Lemma $\ref{processtau}$.
Indeed  for every  $y\geq 0$, $\Psi^{(t)}( \epsilon_{\alpha}(t)(1-\t{m}t)^{-1}y) \ \stackrel{t\rightarrow 1/\t{m}}{\sim} \epsilon_{\alpha}(t)\Psi^{\alpha,1}(y)$. So  ($\ref{equiv}$) entails
$$\epsilon_{\alpha}(t)(1-\t{m}t)^{-1}y \ \stackrel{t\rightarrow 1/\t{m}}{\sim} \ \kappa^{(t)}(-\epsilon_{\alpha}(t)\Psi^{\alpha,1}(y)).$$
Put $y=\kappa^{\alpha,1}(z)$ to get the first limit of the lemma  and follow the same way to get the second one.
\end{proof}
$\newline$

\begin{proof}[Proof of Theorem \ref{asympt}]

First, we prove that $\epsilon_{\alpha}(t).(g(t),d(t))$ converges weakly as $t$ tends to $1/\t{m}$ to $(g(\alpha,1),d(\alpha,1))$.
Indeed by $(\ref{lapdoub})$, we have
$$\E\big(\t{exp}(\rho\epsilon_{\alpha}(t)g(t)-\mu\epsilon_{\alpha}(t)d(t))\big)=(1-\t{m}t)\frac{ \kappa^{(t)}( \epsilon_{\alpha}(t)\rho)-\kappa^{(t)}( \epsilon_{\alpha}(t)\mu)}{ \epsilon_{\alpha}(t)(\rho-\mu)}.$$
Let $t\rightarrow 1/\t{m}$ using Lemma \ref{processtau}  and find the right hand side of (\ref{laplenght}) to conclude.\\ \\

Moreover $\epsilon_{\alpha}(t)\Rpt^{cl}$ (resp. $\epsilon_{\alpha}(t)\Rmt^{cl}$) converges weakly in $\H(\RRR_+)$  as
$t$ tends to $1/\t{m}$ to $\Rpl^{cl}$ (resp. $\Rml^{cl}$). Indeed, by
 Proposition (3.9) in \cite{frm}, this is a consequence of the convergence of the  Laplace exponent of
 $\epsilon_{\alpha}(t)\Rpt$ given by Lemma $\ref{processtau}$. Informally, $\epsilon_{\alpha}(t)\Rpt^{cl}$
 is the range of $\big(\epsilon_{\alpha}(t)\taup^{(t)}_{(1-\t{m}t)\epsilon_{\alpha}(t)^{-1}z} \big)_{z\geq 0}$ whose convergence in $\DD$ follows from
Lemma $\ref{processtau}$.  \\ \\

We can now prove the theorem. We know from  (\ref{reprR}) that
$$ \epsilon_{\alpha}(t) \R(t)=\epsilon_{\alpha}(t).(d(t)+\Rpt) \ \sqcup \ (\w{\epsilon_{\alpha}(t).(-g(t)+\Rmt)})$$
where $\epsilon_{\alpha}(t)\Rmt$,
$\epsilon_{\alpha}(t)(-g(t),d(t))$ and
$\epsilon_{\alpha}(t)\Rpt$ are independent by Theorem
\ref{loi}. Similarly
$$\R(\alpha,1)=(d(\alpha,1)+\Rpl) \ \sqcup \ (\w{-g(\alpha,1)+\Rml})$$
where $\Rml$, $(-g(\alpha,1),d(\alpha,1))$ and $\Rpl$  are independent. As remarked above, we have also the following weak convergences as $t$ tends to $1/\t{m}$ :
$$  \epsilon_{\alpha}(t) \Rmt^{cl} \Rightarrow    \Rml^{cl},
\quad  \epsilon_{\alpha}(t)(-g(t),d(t))\Rightarrow (-g(\alpha,1),d(\alpha,1)), \quad
\epsilon_{\alpha}(t) \Rpt^{cl} \Rightarrow \Rpl^{cl}.$$
So $\epsilon_{\alpha}(t)\R(t)^{cl}$ converges weakly to $\R(\alpha,1)^{cl}$ in $\H(\RRR)$ as $t$ tends to $1/\t{m}$.
\end{proof}
$\newline$

\begin{proof}[Proof of  Corollary \ref{asymptl}]
The first result is a direct consequence of Theorem \ref{asympt}. We have then
$$\epsilon_{\alpha}(t)l(t)\stackrel{t\rightarrow 1/\t{m}}{\Longrightarrow } d(\alpha,1)-g(\alpha,1).$$
Moreover, as $\kappa^{2+,1}\circ (-\Psi^{2+,1})=\rm{Id}$, we can compute $\kappa^{2+,1}$ and $(\ref{laplenght})$ gives
$$ \E\big(\t{exp}(-\mu (d(2+,1))- g(2+,1)) \big)=(\kappa^{2+,1})'(\mu)=
\big(
\frac{1+\sqrt{1+2\frac{\t{m}_2}{\t{m}}\mu}}{\frac{\t{m}_2}{\t{m}}}\big)'(\mu)=\frac{1}{\sqrt{-1+2\frac{\t{m}_2}{\t{m}}\mu}}.$$
So, by identification of Laplace transform,
$d(\alpha,1)-g(\alpha,1)$  is  a gamma variable of parameter
$(1/2,\t{m}/(4\t{m}_2))$ and we get the result. The argument is
similar  in
the case $\alpha=2$.
\end{proof}
$\newline$
\begin{proof}[Proof of Theorem \ref{thtrst}]
The argument is similar to that of the proof  the previous theorem  and use the others limits of Lemma $\ref{processtau}$  to get that
if $x\rightarrow \infty$ and  $1-\t{m}t\sim \lambda f_{\alpha}(x)$ with $\lambda>0$, then   $x^{-1}\R(t)$ converges
weakly in $\H(\RRR)$ to  $\{x \in \RRR: Y^{\alpha,\lambda}_x=I^{\alpha,\lambda}_x\}^{cl}$. The theorem follows by restriction to $[0,1]$.
\end{proof}
$\newline$

To prove the corollary of Theorem \ref{thtrst}, we need the following result.
\begin{Lem}
\label{lemexc}
The largest length of excursion
of $(Y^{\alpha,\lambda}_x-I^{\alpha,\lambda}_x)_{x\in [0,1]}$, denoted
by $B_{\alpha,\lambda}$, converges in probability to $0$ as
$\lambda$ tends to infinity and to $1$ as $\lambda$ tends to $0$.
\end{Lem}
\begin{proof}
$\bullet$ Let $0\leq a<b\leq 1$. Note that for all $\lambda'\geq 1$ and $x\geq 0$, $Y_x^{\alpha,\lambda'}-Y_x^{\alpha,1}
=(1-\lambda')x$ ensures that
$I_x^{\alpha,\lambda'}-I_x^{\alpha,1}\geq (1-\lambda')x$. Then,
$$Y^{\alpha,\lambda'}_{a+2\frac{b-a}{3}}-I^{\alpha,\lambda'}_{a+\frac{b-a}{3}} \leq Y^{\alpha,1}_{a+2\frac{b-a}{3}}
-I^{\alpha,1}_{a+\frac{b-a}{3}}+(1-\lambda')\frac{b-a}{3}.$$
So a.s there exists $\lambda'$ such that
$$Y^{\alpha,\lambda'}_{a+2\frac{b-a}{3}}<I^{\alpha,\lambda'}_{a+\frac{b-a}{3}}.$$
As $Y^{\alpha,\lambda'}$ has no negative jumps, it reaches its
infimum on $]-\infty,2(b-a)/3]$ in a point
$c\in[a+(b-a)/3,a+2(b-a)/3]$. Then  a.s there exists
$c\in[a+(b-a)/3,a+2(b-a)/3]$ and $\lambda'>0$  such that
$c\in\R(\alpha,\lambda')$, which entails that $c$ does not belong
to the interior of $B_{\alpha,\lambda'}$. Adding that
$B_{\alpha,\lambda}$ decreases as $\lambda$ increases, this
property ensures
that $B_{\alpha,\lambda}$ converges in probability to $0$ as $\lambda$ tends to infinity. \\ \\

$\quad $ $\bullet$ As $(Y^{\alpha,0}_x)_{x\in\RRR}$ oscillates when $x$ tends to $-\infty$ (see \cite{lev} Corollary 2 on page 190), then
$$ I_0^{\alpha,\lambda} \stackrel{\lambda\rightarrow 0}{\longrightarrow }-\infty,$$
which ensures that $B_{\alpha,\lambda}$ converges in probability to  $1$ as $\lambda$ tends to $0$.
\end{proof}
$\newline$
\begin{proof}[Proof of  Corollary \ref{cortrst}] The first result is a direct consequence of Theorem \ref{thtrst}. \\

If $o(1-\t{m}t)=f_{\alpha}(x) \ \ (x\rightarrow \infty)$, then for every $\lambda>0$ and $x$ large enough, $t\leq (1-\lambda f_{\alpha}(x))/\t{m}$ and
$$ B_1(x,t)/x\leq  B_1(x,\frac{1-\lambda f_{\alpha}(x)}{\t{m}})/x.$$
The right hand side converges weakly to $B_{\alpha,\lambda}$ as $x$ tends to infinity. Letting $\lambda$ tend to infinity, the lemma above entails that  $B_1(x,t)/x\stackrel{x\rightarrow \infty}{\longrightarrow }0 \ \ \t{in} \ \P$.  \\

Similarly if  \ $1-\t{m}t=o(f_{\alpha}(x)) \ \ (x\rightarrow \infty)$ ,  then for every $\lambda>0$ and $x$ large enough,
$$ B_1(x,t)/x\geq  B_1(x,\frac{1-\lambda f_{\alpha}(x)}{\t{m}})/x.$$
Letting $\lambda$ tend to $0$, Lemma \ref{lemexc} entails  that  $B_1(x,t)/x\stackrel{x\rightarrow \infty}{\longrightarrow }1 \ \ \t{in} \ \P$.
\end{proof}
$\newline$

\noindent{\bf Acknowledgments : } I wish to thank Jean Bertoin for introducing me in this topic and guiding me along the different steps of this work. \\
I am  very grateful to  three anonymous referees  for their insightful comments and for pointing at some important
connections which were missed in the first place.

\end{document}